\documentclass[11pt]{article}
\usepackage{latexsym,amssymb,amsmath,geometry,cite,enumerate,float}
\newtheorem{theorem}{Theorem}[section]
\newtheorem{lemma}[theorem]{Lemma}

\usepackage{setspace}

\begin{document}

\onehalfspace

\title{The Cat and the Noisy Mouse}

\author{Dieter Rautenbach and Moritz Schneider}

\date{}

\maketitle

\begin{center}
Institut f\"{u}r Optimierung und Operations Research, 
Universit\"{a}t Ulm, Ulm, Germany,
\{\texttt{dieter.rautenbach,moritz-1.schneider}\}\texttt{@uni-ulm.de}\\[3mm]
\end{center}

\begin{abstract}
We consider a variant of a pursuit and evasion game studied independently by Britnell and Wildon as well as Haslegrave.
In their game, a cat has to catch an invisible mouse that moves along the edges of some graph $G$. 
In our version, the cat receives partial information about its distance to the mouse, 
and we show that the cat has a winning strategy if and only if $G$ is a forest.
Seager proposed a similar game with complete distance information whose rules cause some small yet important differences to the game we consider.
\end{abstract}

{\small 
\begin{tabular}{lp{13cm}}
{\bf Keywords:} pursuit and evasion game
\end{tabular}
}

\section{Introduction}

\begin{figure}[H]
\begin{center}
\unitlength 0.7mm 
\linethickness{0.4pt}
\ifx\plotpoint\undefined\newsavebox{\plotpoint}\fi 
\begin{picture}(46,39)(0,0)
\put(5,5){\circle*{1.5}}
\put(25,5){\circle*{1.5}}
\put(45,5){\circle*{1.5}}
\put(5,15){\circle*{1.5}}
\put(25,15){\circle*{1.5}}
\put(45,15){\circle*{1.5}}
\put(5,25){\circle*{1.5}}
\put(25,25){\circle*{1.5}}
\put(45,25){\circle*{1.5}}
\put(25,35){\circle*{1.5}}
\put(1,5){\makebox(0,0)[cc]{$u_1$}}
\put(21,5){\makebox(0,0)[cc]{$u_2$}}
\put(49,5){\makebox(0,0)[cc]{$u_3$}}
\put(1,15){\makebox(0,0)[cc]{$v_1$}}
\put(21,15){\makebox(0,0)[cc]{$v_2$}}
\put(49,15){\makebox(0,0)[cc]{$v_3$}}
\put(1,25){\makebox(0,0)[cc]{$w_1$}}
\put(21,25){\makebox(0,0)[cc]{$w_2$}}
\put(49,25){\makebox(0,0)[cc]{$w_3$}}
\put(25,39){\makebox(0,0)[cc]{$x$}}
\put(25,35){\line(0,-1){30}}
\put(45,5){\line(0,1){20}}
\put(45,25){\line(-2,1){20}}
\put(25,35){\line(-2,-1){20}}
\put(5,25){\line(0,-1){20}}
\end{picture}
\end{center}
\vspace{-5mm}\caption{The tree $T^*$.}\label{fig2}
\end{figure}
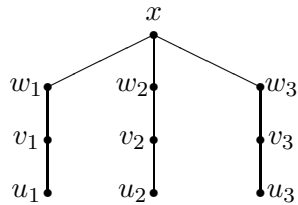

We consider a variant of a pursuit and evasion game independently studied by Britnell and Wildon \cite{brwi} and Haslegrave \cite{ha}. 
The game is played by two players, called the {\it cat} and the {\it mouse} in \cite{ha}, 
on a finite graph $G$ known to both, and proceeds in discrete rounds numbered by positive integers.
In round $i$, the cat chooses a vertex $c_i$ of $G$ and also the mouse chooses a vertex $m_i$ of $G$. 
While the choice of the cat is unrestricted among all vertices of $G$, 
the mouse has to move along an edge of $G$, more precisely, 
the vertices $m_i$ and $m_{i+1}$ have to be distinct yet adjacent for every $i$. 
The cat {\it catches} the mouse in round $i$ if $c_i=m_i$, and the cat wins and terminates the game in the first round in which this happens. 
The mouse wins the game if it can avoid being caught indefinitely, more precisely, 
the mouse wins if for every positive integer $k$, and every sequence $c_1,\ldots,c_k$ of vertices of $G$, 
there is a sequence $m_1,\ldots,m_k$ of vertices of $G$ 
such that $m_{i+1}$ is a neighbor of $m_i$ for every $i\in [k-1]$, and $c_i\not=m_i$ for every $i\in [k]$.
As shown by Britnell and Wildon \cite{brwi} and Haslegrave \cite{ha},
the cat has a winning strategy in their game played on $G$ 
if and only if $G$ is a forest that does not contain the tree $T^*$ shown in Figure \ref{fig2}.

\noindent The key feature of their game is the invisibility of the mouse,
that is, the cat does not know the position of the mouse until it actually catches it. 
In the variant that we consider here, 
the cat receives some information about the distance ${\rm dist}_G(c_i,m_i)$ in $G$ between $c_i$ and $m_i$ in every round $i$;
what information exactly is specified at the beginning of the next section.
Figuratively, the cat estimates its current distance to the mouse using the noise made by the mouse's movement.

Seager \cite{se} introduced a similar game, where the cat, called {\it cop}, knows its exact distance to the mouse, called {\it robber}, 
yet the objective of the mouse is to determine the exact position of the mouse at least once. 
This weaker objective makes the game easier for the mouse; 
as every cop might confirm, knowing where a robber is at some point in time is not quite the same as actually catching him.
There are further differences in Seager's game.
The mouse is forbidden to move to the previous position of the cat but it is allowed not to move at all, that is, 
$m_{i+1}$ is an element of $(\{ m_i\}\cup N_G(m_i))\setminus \{ c_i\}$. 
Allowing $m_{i+1}=m_i$ but not requiring $m_{i+1}\not=c_i$ within Haslegrave's game, 
leads to a trivial situation where any distance information is useless, 
because the mouse can avoid capture indefinitely as soon as the graph has at least one edge. 
Also, requiring $m_{i+1}\in N_G(m_i)\setminus \{ c_i\}$ within Haslegrave's game, 
changes the game considerably; 
if this modified game is played on $T^*$, 
and the cat plays $(c_1,\ldots,c_7)=(x,v_1,x,v_2,x,v_3,x)$,
then, at some point, either the mouse is caught or it has no more legal move. 
In fact, since every second vertex played by the cat is $x$,
in order to avoid capture, 
the mouse must stay in $\{ u_j,v_j,w_j\}$ for some $j\in [3]$. 
Now, when the cat plays $c_i=v_j$, 
the mouse $m_i$ is in $\{ u_j,w_j\}$.
Furthermore, if the mouse $m_i$ is at $w_j$, then it will be caught on $x$ in the next round, and 
if it is at $u_j$, then it has no legal move to play.
It is open what happens to Haslegrave's game if both conditions are imposed, that is, $m_{i+1}\in (\{ m_i\}\cup N_G(m_i))\setminus \{ c_i\}$.

Altogether, Seager's game shows some small yet important differences to the game that we consider here.
As shown in \cite{se,cachdeerwe}, 
the cop wins Seager's game on many graphs that contain cycles, while, in our game, 
the mouse wins as soon as there is any cycle $C$. 
In fact, the mouse's strategy is to stay on $C$ throughout the game, and, 
since every vertex of $C$ has two neighbors on $C$, 
the cat will never be able to guarantee catching the mouse; 
if the cat chooses $c_{i+1}$ as a neighbor of $m_i$, 
then $m_{i+1}$ might just be a neighbor of $m_i$ on $C$ that is distinct from $c_{i+1}$.

Our main result is that the cat has a winning strategy on $G$ exactly if $G$ is a forest.

\section{Result}

We consider Haslegrave's game with partial distance information for the cat.
In what follows, let $d_i$ denote the distance ${\rm dist_G}(c_i,m_i)$ for every positive integer $i$.
In the original game, the cat learns in every round $i$, whether $d_i$ is ``{\it $0$}'' or ``{\it at least $1$}''.
In the version that we consider, the cat receives the following information in every round $i$:
\begin{itemize}
\item whether $d_i$ is ``{\it $0$}'' or ``{\it $1$}'' or ``{\it at least $2$}'', and, additionally,
\item for $i\geq 2$, whether $d_i$ is ``{\it at most $d_{i-1}$}'' or ``{\it bigger than $d_{i-1}$}''.
\end{itemize}
In order to show that the cat has a winning strategy on forests, 
it suffices to describe a winning strategy on trees;
the cat can just apply this strategy to the different components of any forest, one after the other,
and will eventually catch the mouse. 
As we see below, the maximum number of rounds that the mouse can evade capture when the cat plays our strategy on a tree is bounded in terms of the order of the tree. Therefore, if the cat plays on some component of a forest, and this number of rounds is reached without capture, then the mouse must be in another component.

Let $T$ be a tree.

The cat chooses an arbitrary vertex $r$ of $T$ as its {\it root}.
For every vertex $u$ of $T$, let $V_u$ denote the set of vertices of $u$ that are either $u$ or a descendant of $u$.

\medskip

\noindent We consider the game as a sequence of {\it transitions} from some round $i$ to some later round $i+j$. In order to show that the cat has a winning strategy, we argue that until it catches the mouse, it can ensure that all these transitions are of one of four different simple types, which allows the cat to make progress. 

The situation at the beginning of a transition starting with round $i$ is characterized by 
\begin{itemize}
\item a {\it reference vertex} $r_i$, and 
\item two sets $X_i$ and $Y_i$ of vertices.
\end{itemize}
The set $X_i$ contains all vertices that are not a descendant of $r_i$,
which includes $r_i$ itself, 
as well as all vertices in sets $V_u$ for some children $u$ of $r_i$.
Let the children of $r_i$ that do not belong to $X_i$ 
be the leaves $v_1,\ldots,v_k$ and the non-leaves $w_1,\ldots,w_\ell$.
The set $Y_i$ is the union of sets $V_u$ for some children $u$ of $w_1$. Let the children of $w_1$ that do not belong to $Y_i$ be the vertices $x_1,\ldots,x_p$.
Note that $k$, $\ell$, and $p$ are all non-negative integers that might be $0$, and that the cat is free to choose the ordering of the mentioned vertices.

\begin{figure}[H]
\begin{center}
\unitlength 1.5mm 
\linethickness{0.4pt}
\ifx\plotpoint\undefined\newsavebox{\plotpoint}\fi 
\begin{picture}(85,38)(0,0)
\put(34,31){\circle*{1}}
\put(34,31){\line(-3,-1){30}}
\put(34,31){\line(-2,-1){20}}
\put(34,31){\line(0,1){5}}
\put(32,38){\line(1,0){4}}
\put(36,38){\line(0,-1){9}}
\put(36,29){\line(-1,0){3}}
\multiput(33,29)(-.06302521,-.033613445){238}{\line(-1,0){.06302521}}
\put(18,21){\line(0,-1){10}}
\put(0,11){\line(0,1){12}}
\multiput(0,23)(.0978593272,.0336391437){327}{\line(1,0){.0978593272}}
\put(32,34){\line(0,1){4}}
\put(14,21){\circle*{1}}
\put(4,21){\circle*{1}}
\multiput(4,21)(-.03333333,-.11666667){60}{\line(0,-1){.11666667}}
\multiput(31,14)(-.03333333,-.11666667){60}{\line(0,-1){.11666667}}
\multiput(36,14)(-.03333333,-.11666667){60}{\line(0,-1){.11666667}}
\multiput(43,14)(-.03333333,-.11666667){60}{\line(0,-1){.11666667}}
\multiput(53,14)(-.03333333,-.11666667){60}{\line(0,-1){.11666667}}
\multiput(65,14)(-.03333333,-.11666667){60}{\line(0,-1){.11666667}}
\multiput(73,21)(-.03333333,-.11666667){60}{\line(0,-1){.11666667}}
\multiput(83,21)(-.03333333,-.11666667){60}{\line(0,-1){.11666667}}
\put(4,21){\line(0,-1){7}}
\put(31,14){\line(0,-1){7}}
\put(36,14){\line(0,-1){7}}
\put(43,14){\line(0,-1){7}}
\put(53,14){\line(0,-1){7}}
\put(65,14){\line(0,-1){7}}
\put(73,21){\line(0,-1){7}}
\put(83,21){\line(0,-1){7}}
\multiput(4,21)(.03333333,-.11666667){60}{\line(0,-1){.11666667}}
\multiput(31,14)(.03333333,-.11666667){60}{\line(0,-1){.11666667}}
\multiput(36,14)(.03333333,-.11666667){60}{\line(0,-1){.11666667}}
\multiput(43,14)(.03333333,-.11666667){60}{\line(0,-1){.11666667}}
\multiput(53,14)(.03333333,-.11666667){60}{\line(0,-1){.11666667}}
\multiput(65,14)(.03333333,-.11666667){60}{\line(0,-1){.11666667}}
\multiput(73,21)(.03333333,-.11666667){60}{\line(0,-1){.11666667}}
\multiput(83,21)(.03333333,-.11666667){60}{\line(0,-1){.11666667}}
\put(14,21){\line(0,-1){7}}
\multiput(14,21)(.03333333,-.11666667){60}{\line(0,-1){.11666667}}
\multiput(14,21)(-.03333333,-.11666667){60}{\line(0,-1){.11666667}}
\put(24,21){\circle*{1}}
\put(33,21){\circle*{1}}
\put(9,10){\makebox(0,0)[cc]{$X_i$}}
\put(24,18){\makebox(0,0)[cc]{$v_1$}}
\put(33,18){\makebox(0,0)[cc]{$v_k$}}
\put(43,21){\circle*{1}}
\multiput(43,21)(-.057692308,-.033653846){208}{\line(-1,0){.057692308}}
\put(43,21){\line(-1,-1){7}}
\put(31,14){\circle*{1}}
\put(36,14){\circle*{1}}
\put(43,14){\circle*{1}}
\put(53,14){\circle*{1}}
\put(65,14){\circle*{1}}
\put(73,21){\circle*{1}}
\put(83,21){\circle*{1}}
\put(60,14){\makebox(0,0)[cc]{$\ldots$}}
\put(43,14){\line(0,1){7}}
\put(46,23){\makebox(0,0)[cc]{$w_1$}}
\put(39,33){\makebox(0,0)[cc]{$r_i$}}
\put(78,21){\makebox(0,0)[cc]{$\ldots$}}
\multiput(73,21)(-.1313131313,.0336700337){297}{\line(-1,0){.1313131313}}
\put(83,25){\makebox(0,0)[cc]{$w_\ell$}}
\put(73,25){\makebox(0,0)[cc]{$w_2$}}
\put(28,14){\line(0,1){2}}
\put(28,16){\line(1,0){11}}
\put(39,16){\line(0,-1){12}}
\put(28,4){\line(0,1){10}}
\put(34,3){\makebox(0,0)[cc]{$Y_i$}}
\put(39,4){\line(0,-1){4}}
\put(39,0){\line(-1,0){11}}
\put(28,0){\line(0,1){5}}
\multiput(34,31)(.0337078652,-.0374531835){267}{\line(0,-1){.0374531835}}
\put(18,11){\line(0,-1){4}}
\put(18,7){\line(-1,0){18}}
\put(0,7){\line(0,1){5}}
\put(29,21){\makebox(0,0)[cc]{$\ldots$}}
\put(24,21){\line(1,1){10}}
\multiput(34,31)(-.0333333,-.3333333){30}{\line(0,-1){.3333333}}
\multiput(43,21)(.048076923,-.033653846){208}{\line(1,0){.048076923}}
\multiput(83,21)(-.164983165,.0336700337){297}{\line(-1,0){.164983165}}
\put(46,14){\makebox(0,0)[cc]{$x_1$}}
\put(68,14){\makebox(0,0)[cc]{$x_p$}}
\put(56,14){\makebox(0,0)[cc]{$x_2$}}
\multiput(65,14)(-.105769231,.033653846){208}{\line(-1,0){.105769231}}
\end{picture}
\end{center}
\vspace{-5mm}\caption{The situation at the beginning of a transition starting with round $i$.}\label{fig1}
\end{figure}
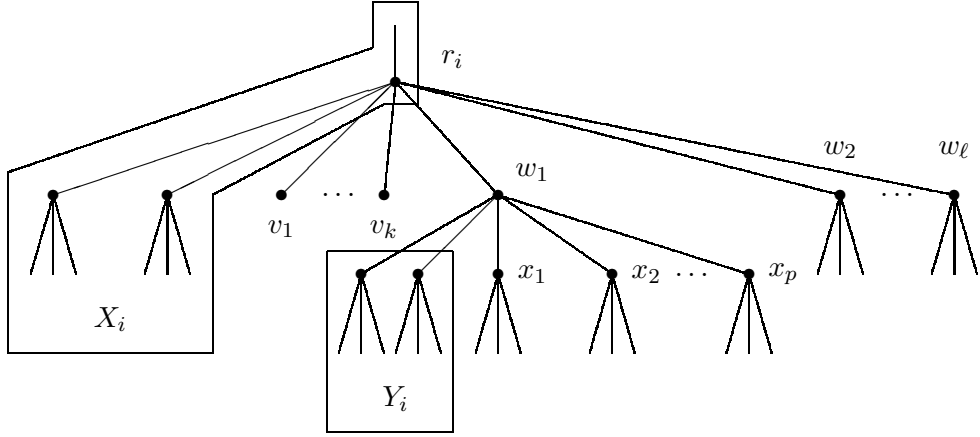

See Figure \ref{fig1} for an illustration.

\noindent In the round $i$ starting the current transition, the cat chooses $c_i$ equal to the reference vertex $r_i$.
Furthermore, exploiting knowledge acquired in earlier rounds,
the cat is sure that $m_i$ does not lie in $X_i\cup Y_i$.
For $i=1$, that is, in the very first round, the reference vertex $r_1$ is the root $r$ of $T$, and the sets $X_1$ and $Y_1$ are both empty,
that is, $m_1\not\in X_1\cup Y_1$ holds trivially.

As said above, we consider four types of transitions from some round $i$ to some later round $i+j$ that are characterized by the following conditions.
\begin{itemize}
\item {\bf Type 1} 

$d_i=1$, 
$j=1$, 
$r_{i+1}$ equals $r_i$, 
$X_i$ is a subset of $X_{i+1}$, 
$Y_{i+1}=\emptyset$, and
$d_{i+1}\geq 2$.
\item {\bf Type 2} 

$d_i\geq 2$, 
$j=1$, 
$r_{i+1}$ is a child of $r_i$, 
$X_i$ is a proper subset of $X_{i+1}$, and 
$Y_{i+1}=\emptyset$.
\item {\bf Type 3} 

$d_i\geq 2$, 
$j=2$, 
$r_{i+2}$ equals $r_i$, 
$X_i$ is a proper subset of $X_{i+2}$, and
$Y_{i+2}=\emptyset$.
\item {\bf Type 4} 

$d_i\geq 2$, 
$j=2$, 
$r_{i+2}$ equals $r_i$, 
$X_i$ is a subset of $X_{i+2}$, 
$Y_i$ is a proper subset of $Y_{i+1}$, and
$d_{i+2}\geq 2$.
\end{itemize}

\begin{lemma}\label{lemma1}
As long as the cat does not win, it can ensure that the game is a sequence of transitions of types 1, 2, 3, and 4.
\end{lemma}
{\it Proof:} Let the situation at the beginning of a transition starting with round $i$ be as summarized above.
As already said, in round $i$ starting the current transition, the cat chooses $c_i$ equal to the reference vertex $r_i$ of that transition, 
and $m_i$ does not lie in $X_i\cup Y_i$. 
Note that the distinction of the different cases discussed below, only requires information that is available to the cat.

Since the cat does not yet win, we have $d_i\geq 1$.

If $d_i=1$, then $m_i\in \{ v_1,\ldots,v_k\}\cup \{ w_1,\ldots,w_\ell\}$.
In this case, the cat chooses $c_{i+1}$ again equal to $r_i$.
Since the mouse must move, either $d_{i+1}=0$ or $d_{i+1}\geq 2$.
Note that this would not hold, if the mouse would be allowed not to move at all like in \cite{se}.
Since we assume that the cat does not yet win, we have $d_{i+1}\geq 2$, which implies $m_{i+1}\not\in X_i$, and, hence, $X_{i+1}$ can be set to $X_i$.
Since $m_i$ might have been $w_1$, and the mouse can move to any child of $w_1$, the set $Y_{i+1}$ is set to empty.
Altogether, setting $j$ to $1$, and not changing the reference vertex $r_{i+1}=r_i$, the cat concludes a transition of type 1.
In what follows, we may assume that $d_i\geq 2$,
which implies that $\ell$ is positive, and that $m_i$ is a descendant of a vertex in $\{ w_1,\ldots,w_\ell\}$ that does not lie in $Y_i$.

If $\ell=1$, then $m_i$ is a descendant of $w_1$.
In this case, the cat chooses $c_{i+1}$ equal to $w_1$.
Now, setting $j$ to $1$,
updating the reference vertex $r_{i+1}$ to the child $w_1$ of $r_i$, 
setting $X_{i+1}$ to $X_i\cup Y_i\cup \{ v_1,\ldots,v_k\}\cup \{ w_1\}$,
and $Y_{i+1}$ to the empty set,
the cat concludes a transition of type 2.
In what follows, we may assume that $\ell\geq 2$.
In this case, the cat chooses $c_{i+1}$ equal to $x_1$.
Since the mouse must move, we have $d_{i+1}\not=d_i$.

If $d_{i+1}\leq d_i$, then $m_i$ was a descendant of $w_1$ outside of $Y_i$. 
Now, setting $j$ to $2$,
not changing the reference vertex $r_{i+2}=r_i$, 
setting $X_{i+2}$ to $X_i\cup \{ v_1,\ldots,v_k\}\cup V_{w_2}\cup \cdots \cup V_{w_\ell}$,
and $Y_{i+2}$ to the empty set,
the cat concludes a transition of type 3.
Note that $m_{i+1}$ might have been $w_1$,
which forces the choice of $Y_{i+2}$.
In what follows, we may assume that $d_{i+1}>d_i$, which implies that
\begin{itemize}
\item either 
$m_i$ is in $V_{w_2}\cup \cdots \cup V_{w_\ell}$,
and $m_{i+1}$ is the parent of $m_i$,
\item or $p\geq 2$, 
$m_i$ is in $V_{x_2}\cup \ldots V_{x_p}$,
and $m_{i+1}$ is a child of $m_i$,
implying ${\rm dist}_G(r_i,m_{i+1})\geq 3$, and, hence, ${\rm dist}_G(r_i,m_{i+2})\geq 2$.
\end{itemize}
Now, the cat sets $j$ to $2$, and does not change the reference vertex $r_{i+2}=r_i$. As specified above, beginning the next transition with round $i+2$, the cat chooses $c_{i+2}$ equal to $r_i$.

If $d_{i+2}=1$, then $m_i$, and hence also $m_{i+1}$, were in $V_{w_2}\cup \cdots \cup V_{w_\ell}$.
Setting $X_{i+2}$ to $X_i\cup \{ v_1,\ldots,v_k\}\cup V_{w_1}$,
and $Y_{i+2}$ to the empty set,
the cat concludes a transition of type 3.
Otherwise, if $d_{i+2}\geq 2$,
then $m_i$ was not a descendant of $x_1$,
setting $X_{i+2}$ to $X_i\cup \{ v_1,\ldots,v_k\}$,
and $Y_{i+2}$ to $Y_i\cup V_{x_1}$,
the cat concludes a transition of type 4.
$\Box$

\begin{theorem}\label{theorem1}
The cat has a winning strategy if and only if the game is played on a forest.
\end{theorem}
{\it Proof:} Since we already explained how the mouse avoids capture if the graph contains a cycle, 
it suffices to show that the cat can win the game played on a tree $T$.
Suppose, invoking Lemma \ref{lemma1}, that the cat ensures that the game is a sequence of transitions of types 1, 2, 3, and 4, 
yet the cat does not win after any finite number of steps.
In order to derive a contradiction, we measure the progress made by the cat using the cardinality of $X_i$, 
where the considered indices $i$ are such that some transition leads from round $i$ to some later round $i+j$.

If the transition is of type 2 or 3, 
then $|X_i|$ increases within that transition, that is, $|X_i|<|X_{i+j}|$.
If the transition is of type $4$, then after at most $n(T)$ further transitions of type 4, there will be a transition of type 2 or 3,
and, again, $|X_i|$ increases.
Finally, every transition of type 1 is followed by a transition of type 2, 3, or 4. Altogether, $|X_i|$ grows after at most every $n(T)$ transitions,
and since it is bounded by $n(T)$, we obtain the desired contradiction.
$\Box$

\medskip

\noindent As a surprising feature of their results, 
Britnell and Wildon \cite{brwi} and Haslegrave \cite{ha} prove the optimality of their strategies for $T^*$-free trees.
We are rather sure that the strategy described in Theorem \ref{theorem1} for our game is not optimal. Nevertheless, the maximum number of rounds that the cat needs to capture the mouse on some tree $T$ using this strategy is at most quadratic in the order $n$ of $T$.

Let $t_i$ denote the number of transitions of type $i$ for $i\in [4]$ during some game.
Clearly, $t_1\leq t_2+t_3+t_4+1$, and $t_2\leq n-1$.
Since for every transition of type 3, the set $X_i$ grows at least by a set of the form $V_w$ for some child $w$ of the current reference vertex $u$, we have 
$t_3\leq \sum\limits_{u\in V(T)}d_T(u)=2n-2.$
Similarly, the number of transitions of type $4$ performed with the same reference vertex $u$ is at most $d_T(u)\left(\sum\limits_{w:\,\,w\,\,is\,\,a\,\,child\,\,of\,\,u}d_T(w)\right),$ and, hence, 
$$t_4\leq \sum_{u\in V(T)}d_T(u)\left(\sum\limits_{w:\,\,w\,\,is\,\,a\,\,child\,\,of\,\,u}d_T(w)\right)
=\sum\limits_{uv\in E(T)}d_T(u)d_T(v)
\leq \left(\sum\limits_{u\in V(T)}d_T(u)\right)^2
=(2n-2)^2.$$
Since every transition of type 1 and 2 consists of one round,
and every transition of type 3 and 4 consists of two rounds, 
the cat wins after at most
$t_1+t_2+2(t_3+t_4)\leq (t_2+t_3+t_4+1)+t_2+2(t_3+t_4)
=2t_2+3t_3+3t_4+1\leq 12n^2-16n+5$
rounds.

\section{Conclusion}

A first obvious question concerns an optimal strategy for our game played on a tree.
Even if the cat knows $d_i$ exactly, 
considering the game played on a path $P$ of order $n$, 
it follows that the mouse can survive about $n/2$ rounds without being caught, 
starting at a vertex close to the middle of $P$, and then, after that, 
moving one step away from the cat, if $c_{i+1}\not=m_i$, and, 
moving to a larger half of $P-m_i$, if $c_{i+1}=m_i$.

Another obvious question is what happens if the cat only learns in round $i$
whether $d_i$ is ``{\it $0$}'' or ``{\it $1$}'' or ``{\it at least $2$}''.
Similarly, it is unclear what happens if the cat only learns
for $i\geq 2$, whether $d_i$ is ``{\it at most $d_{i-1}$}'' or ``{\it bigger than $d_{i-1}$}''.
It is easy to see that for both these variants, the cat has a winning strategy on $T^*$
by choosing, in turn,
$x$ maybe multiple times until it is sure that $d_i\geq 2$ or $d_i>d_{i-1}$, 
then $v_1$,
then again $x$ maybe multiple times until it is sure that $d_{i'}\geq 2$ or $d_{i'}>d_{i'-1}$, 
then $v_2$,
then again $x$ maybe multiple times until it is sure that $d_{i''}\geq 2$ or $d_{i''}>d_{i''-1}$, and
then $v_3$.
Using an initial segment of this strategy, the cat will identify some $j\in [3]$
with $m_i\in \{ u_j,v_j,w_j\}$ for some $i$, and, 
choosing $(c_{i+1},c_{i+2},\ldots)=(x,w_j,v_j,u_j)$ 
necessarily leads to capturing the mouse.

Another open problem is what happens to Haslegrave's game (with or without additional distance information for the cat) 
under the restriction $m_{i+1}\in (\{ m_i\}\cup N_G(m_i))\setminus \{ c_i\}$.

\end{document}